\documentclass{amsart}
\usepackage{algorithm,algpseudocode}
\usepackage[utf8]{inputenc}
\usepackage{amsfonts} 
\usepackage{url}
\usepackage{amsmath} 
\usepackage{amssymb} 
\usepackage{amstext} 
\usepackage{array}   
\usepackage{bbm} 
\usepackage{mathtools} 
\DeclarePairedDelimiter\abs{\lvert}{\rvert} 
\DeclareMathOperator*{\argsup}{arg\,sup}
\DeclareMathOperator{\dom}{dom}
\DeclarePairedDelimiter{\ceil}{\lceil}{\rceil}

\newtheorem{theorem}{Theorem}[section]

\newtheorem{lemma}{Lemma}[section]

\newtheorem{assumption}{Assumption}
\newtheorem{definition}{Definition}[section]
\newtheorem{notation}{Notation}[section]
\newtheorem{property}{Property}[section]

\newtheorem{proposition}{Proposition}[section]
\newtheorem*{remark*}{Remark}
\def\D{\mathcal{D}}
\def\N{\mathbb{N}}
\def\R{\mathbb{R}}
\def\V{\mathbb{V}}

\def\1{\mathbbm1}
\def\lp{\left(}
\def\rp{\right)}
\title[Convergence of D\"umbgen's Algorithm]{Convergence of D\"umbgen's Algorithm for Estimation of Tail Inflation}
\author{Jasha Sommer-Simpson}
\date{\today}
\begin{document}
\begin{abstract}
    Given a density $f$ on the non-negative real line, D\"umbgen's algorithm is a
    routine for finding the (unique) log-convex, non-decreasing function $\hat\phi$ such
    that $\int\hat\phi(x)f(x)dx=1$ and such that the likelihood
    $\prod_{i=1}^{n}f(x_i)\hat\phi(x_i)$ of given data $x_1,\ldots,x_n$ under density
    $x\mapsto \hat\phi(x)f(x)$ is maximized.  We summarize D\"umbgen's algorithm for
    finding this MLE $\hat\phi$, and we present a novel guarantee of the algorithm's
    termination and convergence.
\end{abstract}
\maketitle
\tableofcontents

\section{Introduction}\label{sec:Introduction}

Motivated by the study of statistical sparsity, we suppose that $f$ is a density on the
set $\R_{\geq0}$ of non-negative real numbers, and we consider the class
    \[\Phi_1:=\left\{\phi:\R_{\geq0}\to\R_{>0}~\middle|~\text{$\phi$ log-convex,
    non-decreasing s.t.}\int \phi(x)f(x)dx=1\right\}\]
of positive, log-convex, non-decreasing functions $\phi$ such that $x\mapsto\phi(x)f(x)$
defines a probability density. The interpretation is that, for any given $\phi$ in
$\Phi_1$, the product $\phi\cdot f$ represents the modification of $f$ by the
tail-inflation function $\phi$; the non-decreasing, log-convex nature of $\phi$ means
that the tails of the density $\phi\cdot f$ contain relatively more probability mass
than do the tails of $f$.

Treating $\{\phi\cdot f~:~\phi\in\Phi_1\}$ as a family of densities parameterized by
$\Phi_1$, one defines a maximum likelihood estimate $\hat\phi$ for given data
$x_1,\ldots,x_n\in\R_{\geq0}$ as
\[
    \hat\phi:=\argsup_{\phi\in\Phi_1}\prod_{i=1}^n \phi(x_i)f(x_i).
\]
D\"umbgen~\cite{Dumbgen} has provided an iterative active set algorithm for finding this
maximum-likelihood estimate $\hat\phi$. Despite the fact that the space $\Phi_1$ of
tail inflation functions is infinite-dimensional, D\"umbgen's algorithm is able to
produce a sequence $\phi_0,\phi_1,\ldots,\phi_k,\ldots\in\Phi_1$ of functions that
converges (in likelihood) to the MLE. Making use of the fact that the logarithm
$\hat\theta:=\log(\hat\phi)$ of $\hat\phi$ is piecewise linear and has finitely many
breakpoints (which is proved in Section 5.1 of~\cite{Dumbgen}), Dumbgen's algorithm
iteratively updates an active set $D\subset\R_{\geq0}$ of breakpoints for a convex,
piecewise linear candidate function $\theta$ satisfying $\int e^{\theta(x)}f(x)dx=1$.
After the $k$th iteration of the algorithm we obtain $\phi_k$ by exponentiating the
$k$th candidate function $\theta_k$:
    \[\phi_k(x):=e^{\theta_k(x)}.\]

The main aim of this paper is to establish a guarantee of convergence for
D\"umbgen's algorithm.
The proof presented relies on the following three assumptions:

\begin{assumption}\label{assm:full_support}
    The density $f$ is continuous and has full support on $\R_{>0}$, that is, $f(x)>0$
    for all positive $x$.
\end{assumption}

\begin{assumption}\label{assm:exp_tail}
The density $f$ has an exponential tail, that is,
there exists a constant $\beta\in\R$ such that 
\[
(\forall\lambda\in\R)\quad\quad\quad\lambda<\beta\iff\int e^{\lambda x}f(x)dx<\infty
\]
and
\[
    \lim_{\lambda\to\beta^-}\int e^{\lambda x}f(x)dx = \infty
.
\]
\end{assumption}

\begin{assumption}\label{assm:finite_second_moment}

    Exponential tilting of $f$ results in a density with a finite second moment: for
    $\lambda\in\R$,
    \[
        \int_0^{\infty}e^{\lambda x}f(x)dx<\infty
        \implies
        \int_0^{\infty}x^2e^{\lambda x}f(x)dx<\infty.
    \]
\end{assumption}

The family of Gamma distributions is a prototypical example satisfying the above three
requirements.

It should be noted that D\"umbgen's paper~\cite{Dumbgen} provides an algorithm that
works in the enlarged setting where $f$ is defined on $\R$ and $\phi$ is not required to
be monotone, and that the setting where $\phi$ is log-concave (rather than log-convex)
is also addressed in~\cite{Dumbgen}. The focus of the present paper is restricted to
what Dumbgen calls ``Setting 2B'', where $f$ is defined on $\R_{\geq0}$ and where the
tail-inflation functions $\phi$ are required to be log-convex and non-decreasing. This
being said, the results presented in this paper generalize well to the other settings
considered in~\cite{Dumbgen}.

In Section~\ref{sec:Summary_of_Dumbgen's_paper} we give the statement of D\"umbgen's
algorithm for estimation of $\hat\phi$ in the setting where $f$ is defined on the
non-negative real half-line and $\hat\phi$ is required to be log-convex. Additionally,
we state several key results from D\"umbgen's paper~\cite{Dumbgen} that are important in
demonstrating convergence. 
In section \ref{sec:convergence_proof_overview} we give an overview of the proof of
convergence.
In section \ref{sec:DL_small_implies_theta_eps-optimal} we calculate a bound on the
suboptimality of $\theta_k$ in terms of the directional derivatives of our objective
function $L$ (which is a modified log-likelihood function).

In section \ref{sec:lem_universal_slope_bound} we build on a key result from D\"umbgen's
paper to show that the maximal slope
$\sup_{k,x}\theta_k'(x)$ attained by any candidate function $\theta_k$ is bounded above
by some number strictly less than $\beta$. Additionaly, we give a bound on the values
$\abs{\theta_k(0)}$ taken by the candidate functions at zero.
The results from Section~\ref{sec:lem_universal_slope_bound} are used in section
\ref{sec:bound_on_curvature} to give a lower bound on the change in $L$ resulting from
each step taken by D\"umbgen's algorithm.  In section \ref{sec:proof_of_convergence} we
finish the proof guaranteeing convergence of the algorithm.
In section \ref{sec:conclusion} we conclude.

\section{Statement of Dumbgen's Algorithm and Results from
\cite{Dumbgen}}\label{sec:Summary_of_Dumbgen's_paper}

The aim of this section is to summarize the derivation of D\"umbgen's algorithm for
estimation of a log-convex tail-inflation factor.  We also summarise the results from
D\"umbgen's paper that are used later in proving convergence of the algorithm,
glossing over proofs when convenient.

As mentioned in the introduction, D\"umbgen's paper also considers settings where the
ambient density$f$ is defined on the whole real line (as opposed to on $\R_{\geq0}$),
and where the inflation factor $\phi$ is log-concave instead of log-convex; these
settings are not considered here.  See D\"umbgen's paper for full discussion.

We are given data $x_1,\ldots,x_n\in\R_{\geq0}$ and a density $f$ on $\R_{\geq0}$ with
full support. To find the log-convex, non-decreasing function $\hat\theta$ maximizing
the likelihood of $x_1,\ldots,x_n$, an active-set strategy is
used.

First, we define the set
    $$\Theta_1=\left\{\theta:\R_{\geq0}\to\R~\middle|~\text{$\theta$ convex,
    non-decreasing s.t.}\int e^{\theta(x)}f(x)dx=1\right\}$$
of candidate functions $\theta$. We note the set bijection $\Theta_1\cong\Phi_1$ defined
by $\theta\mapsto e^\theta$.  The log-likelihood of a given candidate
$\theta\in\Theta_1$ is given by
    \begin{equation}\label{phrase:likelihood_of_theta}
    \log\left[\prod_{i=1}^n e^{\theta(x_i)}f(x_i)\right] =\sum_{i=1}^n
    \left[\theta(x_i)+\log f(x_i)\right].
    \end{equation}
Seeing as $f$ and $x_1,\ldots,x_n$ are fixed, we take as our objective for optimiziation
the simplification
    \begin{equation}\label{def:l(theta)}
    l(\theta):=\sum_{i=1}^n \theta(x_i).
    \end{equation}
of the log-likelihood (\ref{phrase:likelihood_of_theta}). Here $l$ defines a function from
$\Theta_1$ to $\R$.

In order to employ techniques from convex optimization, we consider the superset
    $$\Theta=\left\{\theta:\R_{\geq0}\to\R~\middle|~\text{$\theta$ convex,
    non-decreasing}\right\}$$
of $\Theta_1$. This set $\Theta$ is closed under convex combinations, that is if
$\theta_a$ and $\theta_b$ are elements of $\Theta$ then so is
$\lambda\theta_a+(1-\lambda)\theta_b$, so long as $\lambda$ satisfies
$0\leq\lambda\leq1$.

Dumbgen defines the function $L:\Theta\to\bar\R$ by
\begin{equation}\label{def:L(theta)}
    L(\theta):=\sum_{i=1}^n\theta(x_i)-\int e^{\theta(x)}f(x)dx+1
\end{equation}
We note that $L(\theta)$ is finite if and only if $\int e^{\theta(x)}f(x)$ is finite.

Following the notation used in \cite{Dumbgen}, we use $\hat P$ to
denote the empirical distribution $\frac1n\sum_{i=1}^n\delta_{x_n}$ of the observed data
$x_1,\ldots,x_n$, so that the log-likelihood function (\ref{def:l(theta)}) can be written as
\begin{equation*}\label{eqn:l}
    l(\theta)=\int\theta d\hat P.
\end{equation*}
We write $M$ for the measure on $\R_{\geq0}$ having density $f$, so that the
objective (\ref{def:L(theta)}) can be written as
\begin{equation}\label{eqn:L}
    L(\theta)=\int\theta d\hat P-\int e^\theta dM+1.
\end{equation}

The following four properties demonstrate that $L$ is a suitable objective function for
finding the MLE $\hat\theta$.
\begin{property}
        For any $\theta$ in $\Theta$, we have $L(\theta)=\sum_{i=1}^n\theta(x_i)$ if and
        only if $\int e^{\theta(x)}f(x)dx=1$. It follows that we have functional
        equality $l=L|_{\Theta_1}$ between the log-likelihood function $l$ and the
        restriction of $L$ to $\Theta_1$.
    \begin{proof}  
        For any $\theta$ in $\Theta$, we have
        \[\int e^\theta dM=1\iff l(\theta)=\sum_{i=1}^n\theta(x_i)=\sum_{i=1}^n\theta(x_i)-\int e^\theta dM+1=L(\theta)\]
    \end{proof}
\end{property}
\begin{property}
        The function $L$ is strictly concave on the set $\Theta$.
        \begin{proof}
            To see that $L$ is strictly concave, observe that the first term $\int\theta
            d\hat P$ of (\ref{eqn:L}) is affine, and that the exponential function
            appearing in the second term is strictly convex.
        \end{proof}
\end{property}
\begin{lemma}\label{lem:L_has_unique_maximum}
        $L$ attains its maximum at a unique point
            $$\hat\theta:=\argsup_{\theta\in\Theta}L(\theta).$$
        By the previous property, $\hat\theta$ is an element of $\Theta_1$.
\end{lemma}
\noindent For a proof of the above result, see Lemma 2.7 and Section 5.1 from
\cite{Dumbgen}. We remark that D\"umbgen's proof of this fact does not rely on
Assumptions~\ref{assm:exp_tail} or \ref{assm:finite_second_moment}, and that the result
can still be proved even with weakend versions of Assumption \ref{assm:full_support}.
\begin{property}
        The maximum $\hat\theta:=\argsup_{\theta\in\Theta}L(\theta)$ is an element of
        $\Theta_1$.
        \begin{proof}
            This can be seen by letting $c\in\R$ and taking the derivative of
            $L(\theta+c)$ with respect to $c$. We find that $\frac\partial{\partial
            c}L(\theta+c)=0$ only if $c=-\log\int e^\theta dM$ which implies
            $\theta+c\in\Theta_1$.
        \end{proof}
\end{property}

It follows from the above four properties that maximizing $L$ over $\Theta$ is
equivalent to maximizing $l$ over $\Theta_1$, in the sense that
$$\argsup_{\theta\in\Theta}L(\theta)=\hat\theta=\argsup_{\theta\in\Theta_1}l(\theta).$$

D\"umbgen's algorithm relies on the following crucial lemma, which is listed as Lemma
2.7 in D\"umbgen's paper~\cite{Dumbgen} and is proved in section 5.1 of the same.
\begin{lemma}\label{lem:MLE_IS_PIECEWISE_LINEAR}
    The MLE $\hat\theta$ is a piecewise linear function having finitely many
    breakpoints. Writing $x_1\leq\cdots\leq x_n$ without loss of generality, there is at
    most one breakpoint in the open inverval $(x_i,x_{i+1})$ between each pair of
    adjacent observations $x_i,x_{i+1}$. Moreover, every breakpoint of $\hat\theta$ is
    an element of the set $\left(\{0\}\cup[x_1,x_n]\right)\setminus\{x_1,\ldots,x_n\}$,
    where a breakpoint at $0$ is interpreted as saying that the right derivative of
    $\hat\theta$ at $0$ is nonzero. It follows that $\hat\theta$ has at most $n$
    breakpoints.
\end{lemma}

Thus we consider the set
    $$\V:=\left\{v:\R_{\geq0}\to\R~\middle|~\text{$v$ piecewise linear with finitely
    many breakpoints}\right\}$$
of piecewise linear functions with finitely many breakpoints.  In
particular, we know that $\hat\theta$ belongs to the set $\Theta_1\cap\V$.

\begin{notation}
    Given an element $v$ of $\V$, we let $D(v)$ denote the set
        $$D(v):=\{\tau\in\R_{\geq0}~:~v'(\tau-)\neq v'(\tau+)\}$$
    of breakpoints the given function $v$.  In the above display equation, $v'(\tau-)$
    denotes the left derivative of $v$ at $\tau$, and $v'(\tau+)$ denotes the right
    derivative at $\tau$.
\end{notation}
By the definition of $\V$, the set $D(v)$ is a finite subset of $\R_{\geq0}$ for any $v$
in $\V$.

\begin{notation}[The subset $\V_S$ of $\V$]
    Given a finite subset $S$ of $\R_{\geq0}$, we let $\V_S$ denote the subset
        $$V_S:=\{v\in \V~:~D(v)\subseteq S\}$$
    of $\V$ consisting of functions $v$ with breakpoints in $S$.
\end{notation}

\begin{notation}[The function $V_\tau$]
    For a given non-negative real number $\tau$, let $V_\tau:\R_{\geq0}\to\R$ denote the
    function
        $$V_\tau(x):= (x-\tau)^+=\begin{cases}0&x\leq\tau\\x-\tau&\tau\leq x\end{cases}$$
    that is constantly zero on the interval $[0,\tau]$, and that is increasing with unit
    slope on the interval $[\tau,\infty)$.
\end{notation}
We note that, for any non-negative $\tau$, the function $V_\tau$ is an element of the
set $\Theta\cap\V$ of convex non-decreasing functions with finitely many breakpoints.
Also, it is worth mentioning that $ \Theta\cap\V= \text{ span } ^+\{ ( x-\tau ) ^+\} $
is the convex cone consisting precisely of the set of finite linear combinations of
functions $ \{V_\tau~:~\tau\in\R^+\}$, where the coefficients of the terms in the linear
combinations are non-negative.

Dumbgen's algorithm works by maintaining a set $S$ of breakpoints. The algorithm
alternates between a ``Local Search'', which finds
$\hat\theta_S:=\argsup_{\theta\in\Theta\cap\V_S}L(v)$ and replaces $S$ with the subset
$D(\hat\theta_S)$ of $S$, and a ``Global Search'' which replaces the set $S$ with
$S\cup\{\tau\}$, where $\tau$ is chosen as to maximize the directional derivative
    $$\lim_{t\to0^+}\frac{L(\hat\theta_S+tV_\tau)-L(\hat\theta_S)}t$$
in the direction of $V_\tau$. In general, for $\theta$ in
$\Theta$ and $v$ in $\V$, we write
    $$DL(\theta,v):=\lim_{t\to0^+}\frac{L(\theta+tv)-L(\theta)}t.$$

\begin{algorithm}
\caption{D\"umbgen's Algorithm for Finding of Log-Convex Tail Inflation MLE}
\label{alg:Dumbgen}
\begin{algorithmic}[1]
    \Procedure{D\"umbgen'sAlgorithm}{$\theta_0,\delta_0,\delta_1,\epsilon$}
        \State $k\gets 0$
        \Loop
            \State $\tau_k\gets\argsup_\tau DL(\theta_k,V_\tau)$
                \Comment{to within $\delta_0$ suboptimality}
                \label{alg:line:tau_k}
            \State $h_k\gets DL(\theta_k, V_{\tau_k})$
            \If{$h_k\leq\epsilon$}
                \Comment{$h_k\leq\epsilon$ is the termination criterion}
                \label{alg:line:termination_criterion}
                \State \Return $\theta_k$
            \EndIf
            \State $k\gets k+1$
            \State $S_k\gets D(\theta_{k-1})\cup\{\tau_{k-1}\}$
            \Comment{Note that $D(\theta_{k-1})\subseteq S_{k-1}$}
            \State $\theta_k\gets\argsup_{\theta\in\Theta\cap\V_{S_k}}L(\theta)$
                \Comment{to within $\delta_1$ suboptimality}
                \label{alg:line:theta_k}
        \EndLoop
    \EndProcedure
\end{algorithmic}
\end{algorithm}

It is possible to efficiently find $\argsup_\tau DL(\theta_k,V_\tau)$ on line
\ref{alg:line:tau_k} of Algorithm \ref{alg:Dumbgen} because of the following crucial
lemma, proved in Section 3.2 of D\"umbgen's paper \cite{Dumbgen}.
\begin{lemma}
    For any given $\theta\in\Theta$, the function $\tau\mapsto DL(\theta,V_\tau)$ is
    strictly concave on each of the intervals $[x_i,x_{i+1}]$, where $x_1,\ldots,x_n$
    are the observed data sorted in increasing order.
\end{lemma}
\noindent In particular, we can use a concave optimization routine on the intervals
$[x_i,x_{i+1}]$ to find the point $\tau_k$ such that
$DL(\theta_k,\tau_k)\geq\argsup_\tau DL(\theta_k,V_\tau)-\delta_0$.  Thus, if the
termination criterion $h_k\leq\epsilon$ on line \ref{alg:line:termination_criterion} on
Algorithm \ref{alg:Dumbgen} is met, then we can guarantee
    $$\sup_\tau DL(\theta_k,V_\tau)
    \leq DL(\theta_k,\tau_k)+\delta_0
    = h_k+\delta_0
    \leq\epsilon+\delta_0.$$
We prove in Section \ref{sec:DL_small_implies_theta_eps-optimal} that if
    $$\left(\sup_{\theta\in{\Theta\cap\V_{S_{k}}}}
    L(\theta)\right)-L(\theta_k)\leq\delta_1$$
and if
    $$\sup_\tau DL(\theta_k,V_\tau)\leq\epsilon+\delta_0$$
then $L(\hat\theta)-L(\theta_k)\leq\beta(\epsilon+\delta_0)+\delta_1$ for some constant
$\beta$. Thus the termination criterion $h_k\leq\epsilon$ corresponds directly to
suboptimality of the candidate function $\theta_k$ produced by step $k$ of the
algorithm.

The set $\V_{S_k}$ consists of piecewise linear functions having breakpoints in the
set $S_k$.

Finally, we should mention that, because the space $V_{S_k}$ is finite-dimensional
(having dimension $\abs{S_k}$), it is possible to find
$\argsup_{\theta\in\Theta\cap\V_{S_k}}L(\theta)$, as on line \ref{alg:line:theta_k} of
Algorithm \ref{alg:Dumbgen}, using standard convex optimization procedure.  In
D\"umbgen's implementation of the algorithm, Newton's method is used with a
Goldstein-Armijo stepsize correction. In Appendix \ref{bound-on-number-of-breakpoints},
we sketch a proof that the cardinality of $S_k$
is bounded above by $2n-1$.

\subsection{Other useful results from \cite{Dumbgen}}
\begin{lemma}[D\"umbgen]\label{Silverman's-lemma}
    Let $\1:\R_{\geq0}\to\R$ denote the constant function $x\mapsto1$. An element
    $\theta$ of $\Theta$ belongs to $\Theta_1$ if and only if $DL(\theta,\1)=0$, i.e.
        $$\Theta_1=\{\theta\in\Theta~:~DL(\theta,\1)=0\}.$$
    \begin{proof}
        This proof is from Section 1 of \cite{Dumbgen}.
        We consider the derivative $\frac{\partial}{\partial c}L(\theta+c\1)$.
        We have
        \begin{align*}
            DL(\theta,\1)
            &\equiv\frac\partial{\partial t}L(\theta+t\1)\big|_{t=0^+}
          \\&=\frac\partial{\partial t}\left[\int(\theta+t\1)d\hat P -\int
              e^{\theta+t\1}dM+1\right|_{t=0^+}
          \\&=\left[\int\1 d\hat P -\int
              \1 e^{\theta+t\1}dM\right|_{t=0^+}
          \\&=1 -\int
              \1 e^{\theta}dM=1-\int e^\theta dM
        \end{align*}
        which shows that $DL(\theta,\1)=0$ if and only if $1=\int e^\theta dM$, which
        (by definition) holds if and only if $\theta$ is an element of $\Theta_1$.
    \end{proof}
\end{lemma}
\begin{definition}\label{defn:locally_optimal}
    We say that a function $\theta\in\Theta\cap\V$ is \textbf{locally optimal} if
    $\theta$ maximizes $L$ over the set $\Theta\cap\V_{D(\theta)}$, that is, if
        \begin{equation}
        \theta=\argsup_{v\in\Theta\cap\V_{D(\theta)}}L(v).
        \end{equation}
\end{definition}
\begin{lemma}[D\"umbgen]\label{lem:local-optimality-criterion}
    An element $\theta$ of $\Theta\cap\V$ is locally optiomal if and only if $DL(\theta,\1)=0$ and $DL(\theta,V_\tau)=0$ for all $\tau$ in $D(\theta)$, that is,
        $$
        \left(\theta=\argsup_{v\in\Theta\cap\V_{D(\theta)}}L(v)\right)
        \iff
        \biggl(DL(\theta,\1)=0\quad\text{and}\quad DL(\theta,V_\tau)=0~ \forall\tau\in D(\theta)\biggr)
        .$$
    \begin{proof}
        Given an element $\theta$ of $\Theta\cap\V$,
        the set $\V_{D(\theta)}$ is a finite-dimensional vector space with basis $\{\1\}\cup\{V_\tau~:~\tau\in D(\theta)\}$. The set $\Theta\cap\V_{D(\theta)}$ is the convex cone
        $$\Theta\cap\V_{D(\theta)}=\left\{\alpha\1+\sum_{\tau\in D(\theta)}\beta_\tau V_\tau~\middle|~\alpha\in\R,~\beta_\tau\geq0\forall\tau\in D(\theta)\right\}
        .$$
        in the vector space $\V_{D(\theta)}$. We note that since $\Theta\cap\V_{D(\theta)}$ is a convex set, and is a subset of $\Theta$, the restriction of $L$ to $\Theta\cap\V_{D(\theta)}$ is a convex function.
        
        We use the notation $\beta_{\tau,\theta}$ to denote the change in slope of $\theta$ at breakpoint $\tau$, so that $\theta=\theta(0)\1+\sum_{\tau\in D(\theta)}\beta_{\tau,\theta}V_\tau$.
        By the definition of $D(\theta)$, we have $\beta_{\tau,\theta}>0$ for all $\tau\in D(\theta)$. This is to say that $\theta$ lies in the interior of the convex cone $\Theta\cap\V_{D(\theta)}$. Since $\theta$ is not on the boundary of this set, the local optimality of $\theta$ is equivalent to the condition $\frac{\partial L(\theta)}{\partial \1}=0$ and $\frac{\partial L(\theta)}{\partial V_\tau}=0$ for each $\tau$.
    \end{proof}
\end{lemma}

\begin{property}\label{property:theta_locally_opt_implies_theta_in_Theta1}
        If an element $\theta$ of $\Theta\cap\V$ is locally optimal, then $\theta$ is an element of $\Theta_1$.
        \begin{proof}
            This follows directly from Lemmas \ref{Silverman's-lemma} and \ref{lem:local-optimality-criterion}: If $\theta$ is locally optimal then $DL(\theta,\1)=0$, which is equivalent to $\theta$'s membership in the subset $\Theta_1$ of $\Theta$.
        \end{proof}
\end{property}

\section{Proof of Convergence: Overview}\label{sec:convergence_proof_overview}
To simplify our analysis, we suppose that the search steps
    \[
        \tau_k\gets\argsup_\tau DL(\theta_k,V_\tau) \quad\quad\text{and}\quad\quad
        \theta_k\gets\argsup_{\theta\in\Theta\cap\V_{S_{k-1}\cup\{\tau_{k-1}\}}}L(\theta)
    \]
in Algorithm~\ref{alg:Dumbgen} are exact.  As mentioned at the end of the previous
section, if the termination criterion $h_k\leq\epsilon$ from
line~\ref{alg:line:termination_criterion} of Algorithm~\ref{alg:Dumbgen} is satisfied
then
the
bound
    \begin{equation}\label{ineq:bound_suboptimality_of_theta_k}
    L(\hat\theta)-L(\theta_k)<\beta\epsilon
    \end{equation}
follows. This is proved in Section~\ref{sec:DL_small_implies_theta_eps-optimal}.

Now, supposing that the termination criterion $h_k\leq\epsilon$ is not 
met, we show that, for any fixed positive $\epsilon$ and any real number $R$,
if $L(\theta_k)\geq R$ then there exists some constant $C_{R,\epsilon}>0$ such that
    \[
        h_k > \epsilon \quad\implies\quad L(\theta_{k+1})-L(\theta_k)\geq C_{R,\epsilon}
    \]
for each $k$. This will be proved in Section~\ref{sec:bound_on_curvature}, with help
from lemmas proved in Section~\ref{sec:lem_universal_slope_bound}.

Because $\theta_{k}$ is defined as the $\argsup$ of $L$ over the class
$\Theta\cup\V_{S_{k-1}\cup\{\tau_{k-1}\}}$ of convex piecewise linear functions having
breakpoints in the set $S_{k-1}\cup\{\tau_{k-1}\}$, and because $\theta_{k-1}$ belongs
to this same class, we are guaranteed that $L(\theta_k)\geq L(\theta_{k-1})$ for each
$k$.  Therefore $L(\theta_k)\geq L(\theta_0)$ for each $k$, hence
    $$ h_k > \epsilon \quad\implies\quad L(\theta_{k+1})-L(\theta_k)\geq
    C_{L(\theta_0),\epsilon}$$
for each $k$. Therefore, we can guarantee that after finitely many steps (bounded above
in number by the ratio $(L(\hat\theta)-L(\theta_0))/{C_{L(\theta_0),\epsilon}}$) the
bound (\ref{ineq:bound_suboptimality_of_theta_k}) is reached.

\section{$\theta$ locally optimal and $\sup_\tau DL(\theta,V_\tau)\leq\epsilon$ implies
$L(\hat\theta)-L(\theta)<(const)\cdot\epsilon$}\label{sec:DL_small_implies_theta_eps-optimal}

The goal of this section is to show that suboptimality
$L(\hat\theta)-L(\theta_k)<const\cdot\epsilon$ is implied by the termination condition
$\argsup_{\tau}DL(\theta_k,V_\tau)<\epsilon$. In other words, if there is no $\tau$ satisfying
$DL(\theta_k,V_\tau)>\epsilon$ then the suboptimality $L(\hat\theta)-L(\theta_k)$ of
$\theta_k$ must be small.

As mentioned in the previous section, we simplify out analysis by assuming that the
local search step
    $$\theta\gets\argsup_{v\in\Theta\cap\V_{S}}L(v)$$
is exact, that is,
that $\theta_k$ is locally optimal. Not making this assumption, we would instead obtain
a bound $$L(\hat\theta)-L(\theta_k)<\delta_1+(const)\cdot\epsilon$$ on the suboptimality
of $\theta_k$, where $\delta_1$ is the tolerance parameter for local suboptimality of
$\theta_k$:
$$\delta_1\geq \left(\sup_{v\in{\Theta\cap{V_{D(\theta_k)}}}}L(v)\right) - L(\theta_k).$$

We have stated in Assumption~\ref{assm:exp_tail} that the density $f$ has an exponential
tail, that is, there exists some constant $\beta\in\R$ satisfying
    \begin{equation}\label{inequality_beta_1}
        \int_0^\infty e^{\kappa x}M(dx)<\infty \iff \kappa<\beta.
    \end{equation}
Before we can bound the suboptimality of $\theta_k$ directly need a lemma formalizing
the relationship between $\beta$ and the maximal slope obtained by a convex function
$\theta\in\Theta$.
\begin{lemma}
    Suppose that $\theta\in\Theta$. Define $m(\theta):=\sup_x\theta'(x+)$ to be the
    maximal slope attained by $\theta$. We have
    \[
        m(\theta)<\beta\iff\int e^{theta(x)}M(dx)<\infty.
    \]
    \begin{proof}
    Note that
    \[ \theta(x_n)+m(\theta)(x-x_n)\leq\theta(x)\leq\theta(0)+m(\theta)x. \]
        Therefore
    \[
        \int e^\theta dM\leq \int e^{\theta(0)+m(\theta)x}M(dx)=e^{\theta(0)}\int
        e^{m(\theta)x}M(dx)
    \]
    which gives
    \[
        m(\theta)<\beta\implies\int e^\theta dM<\infty,
    \]
    and
    \[
        \int e^\theta dM\geq\int
        e^{\theta(x_n)+m(\theta)(x-x_n)}M(dx)=e^{\theta(x_n)-m(\theta)x_n}\int
        e^{m(\theta)x}M(dx)
    \]
    which gives
    \[
        m(\theta)\geq\beta\implies\int e^\theta dM=\infty.
    \]
    \end{proof}
\end{lemma}

In particular, we note that $\beta$ is is an upper bound for the maximal slope
$m(\hat\theta)=\hat\theta'(x_n)$ attained by $\hat\theta$.  Note that by
Lemma~\ref{lem:MLE_IS_PIECEWISE_LINEAR}, the MLE $\hat\theta$ does not have a breakpoint
at $x_n$, so we may write $\hat\theta'(x_n)$ to refer unambiguously to the derivative
\[\hat\theta'(x_n-)=\hat\theta'(x_n)=\hat\theta'(x_n+)\] of $\hat\theta$ at $x_n$.

We now state and prove the main result of this section.

\begin{lemma}\label{lem_eps_optimality}
    Let $\theta\in\Theta_1$.
    Suppose that $\theta$ is locally optimal, that is,
    $$\theta=\argsup_{v\in\Theta\cap\V_{D(\theta)}}L(\theta).$$
    If $\sup_{\tau}DL(\theta,V_\tau)<\epsilon$  then
        $$L(\hat\theta)-L(\theta)\leq\hat\theta'(x_n)\epsilon<\beta\epsilon.$$
    \begin{proof}
        Define ${v:=\hat\theta-\theta}$. Although $v$ might not be convex, we do have
        ${v\in\V}$, that is, $v$ is a piecewise linear function with finitely many
        breakpoints. For each breakpoint ${\tau\in D(v)}$, let
        ${\beta_{\tau,v}:=v'(\tau+)-v'(\tau-)}$ denote the change in slope of $v$ at
        ${\tau}$. Similarly, we write ${\beta_{\tau,\theta}}$ and
        ${\beta_{\tau,\hat\theta}}$, respectively, for the changes in slope
        ${\beta_{\tau,\theta}:=\theta'(\tau+)-\theta'(\tau-)}$ and
        ${\beta_{\tau,\hat\theta}:=\hat\theta'(\tau+)-\hat\theta'(\tau-)}$ of $\theta$
        and $\hat\theta$.
        \par
        The proof proceeds as follows: first, we shall show that, for any
        $\theta\in\Theta$,
            \begin{equation}\label{first point}
            L(\hat\theta)-L(\theta)\leq DL(\theta,v).
            \end{equation}
        Next, we show that if $\theta\in\Theta_1$ then
            \begin{equation}\label{second point}
            {DL(\theta,v)=\sum_{\tau\in D(v)}\beta_{\tau,v}DL(\theta,V_\tau)}.
            \end{equation}
        Finally, we show that if $\sup_\tau D(\theta,V_\tau)<\epsilon$ then
            \begin{equation}\label{third point}
            \sum_{\tau\in
            D(v)}\beta_{\tau,v}DL(\theta,V_\tau)<\epsilon\hat\theta'(x_n+)<\epsilon\beta.
            \end{equation}
        Combining (\ref{first point}), (\ref{second point}) and (\ref{third point})
        gives the desired result.
        \par
        Validitiy of inequality (\ref{first point}) follows from concavity of $L$ on
        $\Theta$, which gives
            \begin{align*}
            L(\hat\theta)-L(\theta)
            &=\frac{L(\hat\theta)t-L(\theta)t}t
            \\&= \frac{\lp L(\theta)(1-t)+L(\hat\theta)t\rp - L(\theta)}{t}
            \\&\leq \frac{L\lp\theta(1-t)+\hat\theta t\rp-L(\theta)}{t}
            \\&=\frac{L\lp\theta + tv\rp - L(\theta)}t
            \end{align*}
        for all $t$ in the interval $(0,1)$.  Taking the limit as $t$ approaches $0$
        from above, we have
            $$L(\hat\theta)-L(\theta)\leq\lim_{t\to-^+}\frac{L\lp\theta + tv\rp
            - L(\theta)}t\equiv DL(\theta,v)$$
        as required.
        \par
        To demonstrate validity of equation (\ref{second point}), first re-write $v$ as
            $$v =v(0)\1+\sum_{\tau\in D(v)}\beta_{\tau,v}V_\tau$$
        where $\1$ denotes the constant function $(x\mapsto1)$. Note that the changes in
        slope $\beta_{\tau,v}$ can be negative, if
        $\beta_{\tau,\theta}>\beta_{\tau,\hat\theta}$.
        \par
        In Theorem~\ref{thm:linearity_of_DL} from Appendix~\ref{appx:linearity_of_DL} we
        show that the operator $v\mapsto DL(\theta,v)$ is linear (with a caveat
        regarding the domain on which $L$ is finite); using this linearity, we have
            \begin{align*}
            DL(\theta,v)
            &=DL\lp \theta,~v(0)\1+\sum_{\tau\in D(v)}\beta_{\tau,v}V_\tau\rp
            \\&=v(0)DL\lp \theta,~\1\rp
            +\sum_{\tau\in D(v)}\beta_{\tau,v}DL\lp \theta,V_\tau\rp.
            \end{align*}
        Because $\theta$ is assumed to be locally optimal, we have $DL(\theta,\1)=0$ and thus equation (\ref{second point}) follows.
        \par
        Finally, assuming that $\sup_\tau D(\theta,V_\tau)<\epsilon$, and using the fact that
        $$DL(\theta,V_\tau)=0\quad\forall\tau\in D(\theta)$$
        follows from local optimality of $\theta$ (c.f. Lemma \ref{lem:local-optimality-criterion}),
        we have
            \begin{align*}
            \sum_{\tau\in D(v)}\beta_{\tau,v}DL(\theta,V_\tau)
            &=
            \sum_{\tau\in D(\hat\theta)\setminus D(\theta)}\beta_{\tau,v}DL(\theta,V_\tau)+
            \sum_{\tau\in D(\theta)}\beta_{\tau,v}DL(\theta,V_\tau)
            \\&=\sum_{\tau\in D(\hat\theta)\setminus D(\theta)}\beta_{\tau,v}DL(\theta,V_\tau)
            &\left(\text{since}~DL(\theta,V_\tau)=0\quad\forall\tau\in D(\theta)\right)
            \\&=\sum_{\tau\in D(\hat\theta)\setminus D(\theta)}\left[\beta_{\tau,\hat\theta}-\beta_{\tau,\theta}\right]DL(\theta,V_\tau)
            \\&=\sum_{\tau\in D(\hat\theta)\setminus D(\theta)}\beta_{\tau,\hat\theta}DL(\theta,V_\tau)
            &(\text{since}~\tau\notin D(\theta)\implies\beta_{\tau,\theta}=0)
            \\&\leq\sum_{\tau\in D(\hat\theta)}\beta_{\tau,\hat\theta}\cdot\epsilon
            \\&=\epsilon\cdot\hat\theta'(x_n+)<\epsilon\beta.
            \end{align*}
    \end{proof}
\end{lemma}
It follows from the above lemma that, if the termination criterion $h_k\leq\epsilon$ in
Algorithm \ref{alg:Dumbgen} is met, then
    $$L(\hat\theta)-L(\theta_k)\leq\hat\theta'(x_n+)\cdot\epsilon<\beta\cdot\epsilon$$
if we are using exact searches
on lines \ref{alg:line:tau_k} and \ref{alg:line:theta_k} of Algorithm~\ref{alg:Dumbgen},
or
$$L(\hat\theta)-L(\theta_k)\leq\delta_1+\hat\theta'(x_n+)\cdot(\epsilon+\delta_0)<\delta_1+\beta\cdot(\epsilon+\delta_0)$$
when using inexact searches with tolerances $\delta_0$ and $\delta_1$.

\section{Upper Bounds on the Slope and Intercept of $\theta_k$}\label{sec:lem_universal_slope_bound}

In this section we show that, for any constant $R\in\R$, there exist real numbers
numbers $s_R\in(0,\infty)$ and $m_R\in(0,\beta)$ such that, for any $\theta\in\Theta$,
\[R\leq L(\theta)\implies \left(\abs{\theta(0)}\leq
s_R\quad\text{and}\quad\sup_{x}\theta'(x)\leq m_R\right).\]
As in the previous section, we will write
\[
    m(\theta):=\sup_x\theta(x+)
\]
for the maximal slope obtained by a function $\theta$ in $\Theta$.

\begin{lemma}[Dumbgen]\label{lem:bound_on_L(theta)} For any $\theta$ in $\Theta$,
    \[L(\theta)\leq-\log\int e^{\theta(x)-\theta(x_n)}M(dx).\]
    \begin{proof}
        This result is proved in section 5.1 of D\"umbgen's paper~\cite{Dumbgen}. We
        reproduce the proof here: \par Seeing as $\theta$ is non-decreasing on $[0,x_n]$
        and $\hat P$ is the empirical distribution of the data $x_1,\ldots,x_n$, we have
        $\int\theta d\hat P\leq\theta(x_n)$. Thus
        \begin{align}
            \label{ineq:bound_L}
            L(\theta)
            &=\int\theta d\hat P-\int e^\theta dM+1 \\
            &\leq\theta(x_n)-\int e^{\theta(x)+\theta(x_n)-\theta(x_n)}M(dx)+1 \\
            \label{ineq:bound_L_by_theta_x_n}
            &=\theta(x_n)-e^{\theta(x_n)}\int e^{\theta(x)-\theta(x_n)} M(dx)+1 \\
            &\leq\sup_{p\in\R}\left[p-e^p\int e^{\theta(x)-\theta(x_n)} M(dx)\right]+1\\
            &=-\log\int e^{\theta(x)-\theta(x_n)}M(dx).
        \end{align}
    \end{proof}
\end{lemma}
The result below follows directly from the preceding lemma.
\begin{proposition}[Dumbgen]\label{prop:lim_m(theta)-to-infty}
    For any $\theta$ in $\Theta$,
        $$L(\theta)\leq-\log\int e^{m(\theta)\cdot(x-x_n)}M(dx).$$
    \begin{proof}
        Because $\theta$ is convex,
        $
        \theta(x)-\theta(x_n)
        \geq
        m(\theta)\cdot(x-x_n)
        $
        for all $x$.
        \\
        Therefore, $L(\theta)\leq
        -\log\int e^{\theta(x)-\theta(x_n)}M(dx)\leq
        -\log\int e^{m(\theta)\cdot(x-x_n)}M(dx)$.
    \end{proof}
\end{proposition}
We note that
\begin{equation}\label{eqn:lim_m(theta)-to-infty}
    \lim_{m(\theta)\to\infty}-\log\int e^{m(\theta)\cdot(x-x_n)}M(dx)=-\infty.
\end{equation}
Using this fact, we are able to prove the following:
\begin{theorem}\label{thm:L_bdd_below_implies_m_theta_bdd_above}
    For any $R\in\R$ there exists a positive real number $m_R<\beta$
    such that
    \[ R\leq L(\theta)\implies m(\theta)\leq m_R \]
    for all $\theta$ in $\Theta$.
    \begin{proof}
        Fix $R$ in $\R$.
        By (\ref{eqn:lim_m(theta)-to-infty}) there exists a constant
        $m_R$ such that
        \begin{equation}\label{ineq:minus_log_int_e_to_the_m}
        m> m_R\implies-\log\int e^{m\cdot(x-x_n)}M(dx)< R
        \end{equation}
        for all $m$ in $\R$. This is equivalent to the statement that there exists $m_R$
        satisfying

        \begin{equation*}
        R\leq-\log\int e^{m\cdot(x-x_n)}M(dx)\implies m\leq m_R
        \end{equation*}

        for all $m$ in $\R$.
        By Proposition \ref{prop:lim_m(theta)-to-infty}, if $R\leq
        L(\theta)$ then $R\leq-\log\int e^{m(\theta)\cdot(x-x_n)}M(dx)$,
        implying $ R\leq L(\theta)\implies m(\theta)\leq m_R$.

        To guarantee that such a value $m_R$ can be found satisfying $m_R<\beta$,
        we take a closer look at (\ref{ineq:minus_log_int_e_to_the_m}).
        The statement
        \[-\log\int e^{m\cdot(x-x_n)}M(dx)< R\]
        is equivalent to
        \[e^{-R}< e^{-mx_n}\int e^{mx}M(dx).\]
        Therefore, (\ref{ineq:minus_log_int_e_to_the_m}) can be restated as
        \begin{equation}\label{ineq:minus_log_int_e_to_the_m_restated}
            m> m_R
            \implies
            e^{-R}<e^{-mx_n}\int e^{mx}M(dx)
            .
        \end{equation}
        By Assumption \ref{assm:exp_tail} we have
        \[
            \lim_{\lambda\to\beta^-}e^{-\lambda x_n}\int e^{\lambda x}f(x)dx = \infty,
        \]
        which goes to show that, for any $R$, we can find a number $m_R<\beta$ such that
        (\ref{ineq:minus_log_int_e_to_the_m_restated}) is satisfied.
    \end{proof}
\end{theorem}
We now go through a similar argument to obtain a bound on $\abs{\theta(0)}$.
\begin{theorem}[Dumbgen]\label{thm:L_bdd_below_implies_theta_0_bdd_above}
    For any $R\in\R$, there exists a non-negative real number $s_R$
    such that
    \[ R\leq L(\theta)\implies \abs{\theta(0)}\leq s_R \]
    for all $\theta$ in $\Theta$.
    \begin{proof}
        We will first obtain a bound on $\theta(x_n)$ of the form
        \[ R\leq L(\theta)\implies \abs{\theta(x_n)}\leq q_R . \]
        We will then combine the bounds $q_R$ and $m_R$ to obtain a bound $s_R$ on
        $\theta(0)$.

        To begin, fix $R$ in $\R$.
        From~(\ref{ineq:bound_L_by_theta_x_n}) we have
        \begin{equation}\label{ineq:bound_L_by_theta_x_n_ref}
            L(\theta) \leq\theta(x_n)-e^{\theta(x_n)}\int_0^\infty
            e^{\theta(x)-\theta(x_n)}M(dx)+1.
        \end{equation}
        Since $\theta(x)-\theta(x_n)\geq0$ for $x$ larger than $x_n$, we have
        \[
            \int_0^\infty e^{\theta(x)-\theta(x_n)}M(dx)\geq\int_{x_n}^\infty M(dx)
        \]
        which, combined with (\ref{ineq:bound_L_by_theta_x_n_ref}), gives the following
        inequality:
        \[
            L(\theta)\leq\theta(x_n)-e^{\theta(x_n)}\int_{x_n}^\infty M(dx)+1.
        \]
        We note that
        \begin{equation}
            \lim_{\abs{\theta(x_n)}\to\infty}
            L(\theta)\leq\theta(x_n)-e^{\theta(x_n)}\int_{x_n}^\infty M(dx)+1
            =-\infty
        \end{equation}
        so that, by the same argument as in the proof of
        Theorem~\ref{thm:L_bdd_below_implies_m_theta_bdd_above},
        there exists a constant $q_R$ satisfying
        \[ R\leq L(\theta)\implies \abs{\theta(x_n)}\leq q_R \]
        for all $\theta$ in $\Theta$.

        Suppose now that $R\leq L(\theta)$, so that $\abs{\theta(x_n)}\leq q_R$ and
        $m(\theta)\leq m_R$.
        Because $\theta$ is non-decreasing, we have
        $\theta(0)\leq\theta(x_n)\leq q_R$. Because $\theta$ is convex, we have
        $-q_R-m_R\cdot x_n\leq\theta(x_n)-m_R\cdot x_n\leq\theta(0)$.
        Combining these two inequalities
        \[ -q_R-m_R\cdot x_n\leq\theta(0)\leq q_R\]
        we obtain
        \[ \abs{\theta(0)}\leq q_R+m_R\cdot x_n.\]
        Thus the required bound on $\abs{\theta(0)}$ is given by
        $s_R:=q_R+m_R\cdot x_n$.

    \end{proof}
\end{theorem}

Recall that in Lemma \ref{lem:L_has_unique_maximum} we cited D\"umbgen's proof that
there exists a unique maximizer $\hat\theta$ for the function $L$.
Although we do not attempt to prove this in the present paper, we remark that the bounds
$m_R$ and $s_R$ derived in this section can be used given an upper bound for $L$:

\begin{remark*}
    For any $ R\in\R $ there is a constant $ U_R $ such that
    \[ R\leq L (  \theta  )  \implies L (  \theta  ) \leq U_R.  \]
    It follows that for $\theta\in\Theta$ and for any $ R\in\R$
    \[ L (  \theta  ) \leq \max ( R,U_R ) \]
    so that the objective function $ L$ is bounded above.
    \begin{proof}
        Suppose $ R\leq L (  \theta  )  $ . Then
        \begin{align*}
            L(\theta)&\leq-\log\int_{0}^\infty e^{m(\theta)(x-x_n)}M(dx)
            \\&=-\log\left[\int_{0}^{x_n}e^{m(\theta)(x-x_n)}M(dx)+\int_{x_n}^{\infty}e^{m(\theta)(x-x_n)}M(dx)\right]
            \\&=-\log\left[\int_{0}^{x_n}e^{m_R(x-x_n)}M(dx)+\int_{x_n}^{\infty}e^{0(x-x_n)}M(dx)\right]
            \\&=const_R.
        \end{align*}
    \end{proof}
\end{remark*}

\section{A local bound on the second directional derivative of $L$}
\label{sec:bound_on_curvature}

In this section we show that, supposing $DL(\theta_{k-1},V_{\tau_{k-1}})\geq\epsilon$,
the improvement $L(\theta_{k-1})-L(\theta_{k})$ at step $k$ is bounded below by a
constant.  Define $R:=L(\theta_{k-1})$ so that by
Theorems~\ref{thm:L_bdd_below_implies_m_theta_bdd_above} and
\ref{thm:L_bdd_below_implies_theta_0_bdd_above} there are numbers $s_R$ and $m_R$
satisfying
\[
    \abs{\theta_{k-1}(0)}\leq s_R
    \quad\text{and}\quad
    \sup_x\theta_{k-1}'(x+)\leq m_R<\beta.
\]
%
Note that
\begin{equation}\label{eqn:lower_bd_for_theta_k-theta_k-1}
    L(\theta_k)-L(\theta_{k-1})\geq L(\theta_{k-1}+tV_{\tau_{k-1}})-L(\theta_{k-1})
\end{equation}
for every real number $t$. This follows from the fact that $\theta_k$ is defined as
\[
    \theta_k \equiv
    \argsup_{\theta\in\Theta\cap\V_{D(\theta_{k-1})\cup\{\tau_{k-1}\}}} L(\theta)
\]

and that, for each $t$, the function $\theta_{k-1}+tV_{\tau_{k-1}}$ is also a member of
the set $\Theta\cap\V_{D(\theta_{k-1})\cup\{\tau_{k-1}\}}$.

Given inequality (\ref{eqn:lower_bd_for_theta_k-theta_k-1}), we can bound
$L(\theta_k)-L(\theta_{k-1})$ below by finding a lower bound for
$\sup_{t\geq0}L(\theta_{k-1}+tV_{\tau_{k-1}})-L(\theta_{k-1})$.  For
notational convenience, we define the function
    $$g(t):=L(\theta_{k-1}+tV_{\tau_{k-1}})-L(\theta_{k-1}).$$

    To find a lower bound for $\sup_{t\geq0}g(t)$, we first note that
    $$g'(t)=\frac\partial{\partial t}\left[L(\theta_{k-1}+tV_{\tau_{k-1}})\right]$$
so that $g'(0)\equiv DL(\theta_{k-1},V_{\tau_{k-1}})$.  Moreover, by strict concavity
of $L$, our function $g$ is strictly concave in $t$, that is, $g''(t)<0$ so long as $t$
satisfies $g(t)>-\infty$.

Below, assuming $g'(0)=DL(\theta_{k-1},V_{\tau_{k-1}})>\epsilon$, we argue that there
exists $T>0$ such that $g(T)\leq0$. This, together with strict concavity of $g$ and the
facts $g(0)=0$ and $g'(0)>0$, go to show that $g$ obtains its supremum at a unique point
$t^*\in(0,T)$.

The details of this argument rely on the fact that the maximal slope operator
$m:\Theta\cap\V\to\R$ is additive: for any $t\geq0$ we have
\[ m(\theta_{k-1}+tV_{\tau_{k-1}})=m(\theta_{k-1})+t\cdot
m(V_{\tau_{k-1}})=m(\theta_{k-1})+t.  \]
Defining $R=L(\theta_{k-1})$, we have by
Theorem~\ref{thm:L_bdd_below_implies_m_theta_bdd_above} that there exists a constant
$m_R$ such that
\[
    m(\theta_{k-1})+t>m_R\implies L(\theta_{k-1}+tV_{\tau_{k-1}})<L(\theta_{k-1}).
\]
Writing $T:=m_R-m(\theta_{k-1})$
gives the desired property $g(T)=L(\theta_{k-1}+TV_{\tau_{k-1}})-L(\theta_{k-1})\leq0$,
and so the supremum
\[t^*:=\argsup_tg(t)\]
must exist.
The following lemma allows us to find a quadratic function $y(t)$ that bounds $g(t)$
below on the set $[0,T]=\{t~:~g(t)\geq0\}$, giving us a lower bound on $g(t^*)$.

\begin{lemma}\label{lem:bound_on_g_prime_prime}
    Suppose $\theta\in\Theta$ and $\tau\in\R_{\geq0}$ satisfy
    $DL(\theta,V_\tau)>\epsilon$, and that $R$ is a parameter satisfying $R\leq
    L(\theta)$ so that
    $\abs{\theta(0)}\leq s_R\quad\text{and}\quad\sup_x\theta'(x+)\leq m_R<\beta.$
    Let $g:\R^+\to\R$ denote the function $t\mapsto L(\theta+tV_\tau)-L(\theta)$.
    For $t$ in the set $\{t~:~g(t)\geq0\}$, the magnitude $\abs{g''(t)}$ of the second
    derivative of $g''$ is bounded above by a constant depending only on $s_R$ and
    $m_R$.
    \begin{proof}
        If $t$ satisfies $0\leq g(t)$ then $R\leq L(\theta+tV_{\tau})$ because
        \[
            0\leq g(t)=L(\theta+tV_\tau)-L(\theta)=L(\theta+tV_\tau)-R.
        \]
         Therefore, for all such $t$ we have 
            $$\theta(0)+tV_{\tau}(0)\leq
            s_R\quad\text{and}\quad\sup_x\left(\theta'(x+)+tV_{\tau}'(x+)\right)\leq
            m_R.$$
        This gives $\theta(x)+tV_\tau(x)\leq s_R+m_Rx$ for all $x$, so that
\begin{align}
\nonumber
\abs{g''(t)}
&=\left\lvert\frac{\partial^2}{\partial t^2} \left[L(\theta+tV_\tau)-L(\theta)\right]\right\rvert
=\left\lvert\frac{\partial^2}{\partial t^2} \left[L(\theta+tV_\tau)\right]\right\rvert
\\
\nonumber
&=\left\lvert\frac{\partial^2}{\partial t^2} \left[\int\left(\theta+tV_\tau\right)d\hat P -\int e^{\theta+tV_\tau}dM+1\right]\right\rvert
\\
\nonumber
&=\left\lvert\frac{\partial}{\partial t} \left[\int V_\tau d\hat P -\int V_\tau e^{\theta+tV_\tau}dM\right]\right\rvert
\\
\nonumber
&=\left\lvert-\int V_\tau^2 e^{\theta+tV_\tau}dM\right\rvert
=\int_0^\infty V_\tau^2 e^{\theta+tV_\tau}dM
\\
\nonumber
&\leq\int_0^\infty V_\tau^2 e^{s_R+m_Rx}dM
\\
\nonumber
&=\int_\tau^\infty (x-\tau)^2 e^{s_R+m_Rx}M(dx)
\\
\nonumber
&=\int_\tau^\infty (x^2-2x\tau+\tau^2) e^{s_R+m_Rx}M(dx)
\\
\nonumber
&\leq\int_\tau^\infty (x^2+\tau^2) e^{s_R+m_Rx}M(dx)
\\
\nonumber
&\leq\int_\tau^\infty (x^2+x_n^2) e^{s_R+m_Rx}M(dx)
\\
\nonumber
&\leq\int_0^\infty (x^2+x_n^2) e^{s_R+m_Rx}M(dx)
\\
\nonumber
&\leq\int_0^\infty (x^2+x_n^2) e^{s_R+m_Rx}M(dx)
\\
&= e^{s_R}\int_0^\infty x^2e^{m_Rx}M(dx)+x_n^2e^{s_R}\int_0^\infty e^{m_Rx}M(dx)
\label{ineq:abs_g_prime_prime}
\end{align}
            for all $t$ satisfying the hypothesis $g(t)\geq0$.

            Since $m_R<\beta$, Assumption \ref{assm:finite_second_moment} gives that the
            last line (\ref{ineq:abs_g_prime_prime}) above is finite.
    \end{proof}
\end{lemma}
By the lemma above, $L(\theta_k)-L(\theta_{k-1})$ is greater than or equal to the
supremum attained by the parabola $y(t)=\epsilon x-\gamma x^2/2$, where
$\gamma=\int_0^\infty (x^2+x_n^2)e^{s_R+m_Rx}M(dx)$.
The peak of the parabola $y$ is attained at $t=\frac\epsilon\gamma$, and so
$\sup_{t}y(t)=\frac{\epsilon^2}{2\gamma}$. Thus we have:

\begin{theorem}
\label{thm:L_theta_k_minus_L_theta_k_minus_1_bdd}
$L(\theta_k)-L(\theta_{k-1})\geq\frac{\epsilon^2}{2\gamma}$
if $ \argsup_\tau DL (  \theta_{k-1},V_\tau ) \geq \epsilon $. 
\end{theorem}
We emphasize that this
bound holds for all $k$ such that the termination criterion
$\sup_\tau DL(\theta_k,V_\tau)<\epsilon$ for D\"umbgen's algorithm has not been met.

\section{Proof of Convergence}
\label{sec:proof_of_convergence}

We produce a bound $K_\epsilon\in\N$ such that, starting with a guess
$\theta_0\in\Theta_1\cap\V$, D\"umbgen's Algorithm is guaranteed to converge in
likelihood within $K_\epsilon$ steps:
\[ L(\hat\theta)-L(\theta_{K_\epsilon})<\beta\cdot\epsilon.  \]
In the following theorem, we let $R:=L(\theta_0)$ denote the objective value of
$\theta_0$ and let $s_R\in(0,\infty)$ and $m_R\in(0,\beta)$ denote constants satisfying
\[ \abs{\theta_k(0)}\leq s_R\quad\text{and}\quad\sup_x\theta_k'(x)\leq m_R \]
for all $k$ in $\N$. Such bounds $m_R$ and $s_R$ are guaranteed to exist by
Theorems~\ref{thm:L_bdd_below_implies_m_theta_bdd_above} and
\ref{thm:L_bdd_below_implies_theta_0_bdd_above}, respectively.

We let $h_0$ denote the maximal directional derivative
\[ h_0:=\sup_\tau DL(\theta_0,V_\tau) \],
and we assume that the search procedures on lines~\ref{alg:line:tau_k} and
\ref{alg:line:theta_k} of Algorithm~\ref{alg:Dumbgen} are exact.

\begin{theorem}
    Let $\gamma_R$ denote the constant
    \[ \gamma_R:= \int_0^\infty (x^2+x_n^2) e^{s_R+m_Rx}M(dx) \]
    derived in Lemma~\ref{lem:bound_on_g_prime_prime}. D\"umbgen's Algorithm reaches
    suboptimality within $K_\epsilon$ steps,
    \[ L(\hat\theta)-L(\theta_{K_\epsilon})<\beta\cdot\epsilon,  \]
    where
    $K_\epsilon
    =\ceil*{\frac{L(\hat\theta)-L(\theta_0)-\beta\epsilon}{\epsilon^2/(2\gamma_R)}}
    \leq\ceil*{\frac{\beta(h_0-\epsilon)}{\epsilon^2/(2\gamma_R)}}
    $.
    \begin{proof}
        By Lemma~\ref{lem_eps_optimality}, if there is $k$ such that
        \begin{equation}\label{ref_subopt_criterion}
        \argsup_\tau DL(\theta_k,V_\tau)\leq \epsilon
        \end{equation}
        then $L(\hat\theta)-L(\theta_k)\leq\beta\epsilon$ as required.

        Defining $K_\epsilon$ as in the statement of the Theorem above, suppose that
        the suboptimality criterion~(\ref{ref_subopt_criterion}) has not been met for
        any of the first $K_\epsilon$ steps taken by the algorithm, that is, suppose
        \[
            \argsup_\tau DL(\theta,V_\tau)\geq\epsilon\quad\quad\quad(\forall
            k<K_\epsilon)
            .
        \]
        Then by Theorem~\ref{thm:L_theta_k_minus_L_theta_k_minus_1_bdd},
        \[
        L(\theta_{k+1})-L(\theta_{k})\geq\frac{\epsilon^2}{2\gamma_R}
            \quad\quad\quad(\forall k<K_\epsilon)
        .\]
        Thus we may derive
        \begin{align*}
            L(\hat\theta)-L(\theta_{K_\epsilon})
            &=L(\hat\theta)-L(\theta_0)-\sum_{k=1}^{K_\epsilon}\bigl[L(\theta_k)-L(\theta_{k-1})\bigr]\\
            &\leq
            L(\hat\theta)-L(\theta_0)-\sum_{k=1}^{K_\epsilon}\frac{\epsilon^2}{2\gamma_R}\\
            &=
            L(\hat\theta)-L(\theta_0)-\ceil*{\frac{L(\hat\theta)-L(\theta_0)-\beta\epsilon}{\epsilon^2/(2\gamma_R)}}\frac{\epsilon^2}{2\gamma_R}\\
            &\leq
            \beta\cdot\epsilon
        \end{align*}
        to complete the proof.
    \end{proof}
\end{theorem}
We note here that in practice, the algorithm appears to converge quite quickly.

\section{Conclusion}
\label{sec:conclusion}

A key result proved by D\"umbgen in deriving this algorithm is that the logarithm
$$\hat\theta(x):=\log(\hat\phi(x))$$ is necessary piecewise linear, with at most one
breakpoint between each pair $x_i,x_j\in\{x_1,\ldots,x_n\}$ of adjacent observations.
This motivates the decision to optimize over the space $\Theta$ of convex functions with
finitely many breakpoints.

This space $\Theta$ is infinite-dimensional.  Despite this, we have been  able to prove
that D\"umbgen's algorithm produces a sequence
$\theta_0,\ldots,\theta_k,\theta_{k+1},\ldots$
of functions that converges to the optimum $\hat\theta$.  The proof of
convergence is made possible by the following characteristics of the problem:
    \begin{itemize}
        \item The optimum $\hat\theta$ is piecewise linear and has finitely many
            breakpoints. For any given set $D$ of breakpoints, we can use a
            finite-dimensional convex optimization routine to find
            $\argsup_{v\in\V_D}L(v)$. Thus, this problem of finding $\hat\theta$ lends
            itself to an active-set approach.
        \item It is possible to efficiently find a good candidate $\tau\in\D$ for
            addition to the active set of breakpoints.
        \item Given $\epsilon>0$, if $L(\hat\theta)-L(\theta_k)\geq\epsilon\beta$ then
            the directional derivative ${\frac\partial{\partial
            t}L(\theta+tV_\tau)\rvert_{t=0+}}$ must be larger than $\epsilon$ (c.f.
            Lemma \ref{lem_eps_optimality}).
        \item Strict concavity of $L$ guarantees that the second directional derivative
            ${\frac{\partial^2}{\partial t^2}L(\theta+tV_\tau)}$ is negative for all
            $t$. Although this second derivative is unbounded below, restriction to the
            level set $\{t ~:~ L(\theta+tV_\tau)>R\}$, where $R$ is an arbitrary real
            number, will allow us to produce a bound on ${\frac{\partial^2}{\partial
            t^2}L(\theta+tV_\tau)}$ that is uniform for different values of $\theta$ and
            $\tau$.  This is be the key to producing a lower bound on the improvement
            $L(\theta_{k+1})-L(\theta_k)$ in objective value (as in
            Theorem~\ref{thm:L_theta_k_minus_L_theta_k_minus_1_bdd}).
    \end{itemize}

\appendix

\section{Sketch of proof that each $\theta_k$ has at most $2n-1$ breakpoints}
\label{bound-on-number-of-breakpoints}

We will show that the locally-optimal parameter $\theta$ returned by the
\text{LocalSearch} procedure as defined in \cite{Dumbgen} will have at most $2n$
breakpoints. It will follow that, since the global parameter search results in the
addition of just one breakpoint, candidate functions considered by the algorithm are
limited to at most $2n+1$ breakpoints.

In D\"umbgen's paper it is also proved that the optimal $\hat\theta$ has at most $n$
breakpoints.

The proof of the claim that \text{LocalSearch} results in at most $2n$ breakpoints is a
consequence of the fact that a locally-optimal parameter can have at most $2$
breakpoints between any two adjacent observations. Formally, we have the following

\begin{theorem}[Boundedness of number of breakpoints]
    Let $x_i,x_j\in\{x_1,\ldots,x_n\}$ be adjacent observations, so that $x_i<x_j$ and
    $\nexists x'\in\{x_1,\ldots,x_n\}$ such that $x_i<x'<x_j$.  For any element
    $\theta\in\V\cap\Theta$, the result
        $$\theta':=\arg\max_{v\in\V_{D(\theta)}\cap\Theta_1}L(v)$$
    of a local search over the set $\V_{D(\theta)}\cap\Theta_1$ has at most two
    breakpoints between $x_i$ and $x_j$, i.e.
        $$\left\lvert D(\theta')\cap(x_i,x_j)\right\rvert\leq2.$$
    \begin{proof}
        Suppose that $\theta$ has more than two breakpoints between $x_i$ an $x_j$, i.e.
        $$\left\lvert D(\theta)\cap(x_i,x_j)\right\rvert\geq3.$$ It will suffice to show
        that $\theta$ is not locally optimal, i.e. that there exists some element
        $v\in\V_{D(theta)}$ satisfying
            \begin{equation}\label{optimality-condition}
            DL(\theta,v)>0.    
            \end{equation}
        Let $\tau_1,\tau_2,\tau_3\in D(\theta)\cap(x_i,x_j)$ be three distinct
        breakpoints of $\theta$ on the interval $(x_i,x_j)$. Without loss of generality
        we suppose that $\tau_1<\tau_2<\tau_3$. Of course, we must have
        $\beta_{\tau_1},\beta_{\tau_2},\beta_{\tau_3}>0$.
        
        We claim that there is an element $v$ of $\V_{\{\tau_1,\tau_2,\tau_3\}}$
        satisfying the condition (\ref{optimality-condition}) above.  Indeed, we can
        define \def\tott{\tau_1,\tau_2,\tau_3} $V_{\tott}\in\V{\{\tott\}}$ by
            $$V_{\tott}(x):= \begin{cases} 
                              0 & x\leq \tau_1 \\
                              -\frac{x-\tau_1}{\tau_2-\tau_1} & \tau_1\leq x\leq\tau_2\\
                              -\frac{\tau_3-x}{\tau_3-\tau_2} & \tau_2\leq x\leq\tau_3\\
                              0 & x\geq \tau_3
                             \end{cases}.  $$
        Let $\gamma$ be any positive number. We find that because $V_{\tott}(x_k)=0$ for
        all observations $x_k\in\{x_1,\ldots,x_n\}$, we have $\int\theta d\hat
        P=\int(\theta+\gamma V_{\tott})d\hat P$.
        Moreover,
            $$\theta(x)+\gamma V_{\tott}(x)<\theta(x)~~~\forall x\in(x_i,x_j)$$
        because $V_{\tott}$ is negative everywhere on $(x_i,x_j)$.  Therefore, we have
            $$\int e^{\theta+\gamma V_{\tott}}dM\leq\int e^\theta dM.$$ We conclude that
        (\ref{optimality-condition}) does indeed hold for $v=V_{\tott}$.
    \end{proof}
\end{theorem}

\section{Linearity of the operator $DL(\theta,-)$}\label{appx:linearity_of_DL}

The goal of this section is to show that the map $v\mapsto DL(\theta,v)$ is linear,
provided that $L$ is finite. First, we state a Lemma.

\begin{lemma}\label{lem:Rockafellar}
    Suppose that $y:\R\to\bar\R$ is a concave function and that $r\in\dom y$, that is,
    $y(r)$ is finite. Suppose that $s\in\R$ and that there exists $\epsilon>0$ such that
    $y(r+\epsilon s)$ is finite.
    Then the one-sided derivative
    \[
        \frac\partial{\partial
        t}y(r+ts)\Bigr|_{t=0^+}\equiv\lim_{t\downarrow0}\frac{y(r+ts)-y(r)}t
    \]
    exists and is finite.
    \begin{proof}
        See the proof of Theorem 23.1 from Rockafellar's book~\cite{Rockafellar}.
    \end{proof}
\end{lemma}

Now we can state the main result of this appendix:

\begin{theorem}\label{thm:linearity_of_DL}
    Take $\theta\in\Theta$ such that $L(\theta)>-\infty$. Let $v_1,\ldots,v_K\in\V$ such
    that there exists $\epsilon>0$ satisfying
    \begin{equation}\label{condition:v_k_direction_in_domain}
        \theta+\epsilon v_k\in\Theta
        \quad
        \bigl(\forall k\in\{1,\ldots,K\}\bigr)
        .
    \end{equation}
    Then for any non-negative real coefficients $h_1,\ldots,h_K\geq0$ we have:
    \begin{enumerate}
        \item $DL(\theta,v_k)$ exists and is finite for each $k$ in $\{1,\ldots,K\}$,
            \label{appxstmt:1}
        \item $DL(\theta,\sum_{k=1}^K h_kv_k)$ exists and is finite, and
            \label{appxstmt:2}
        \item there is equality $\sum_{k=1}^K h_k DL(\theta,v_k)=DL(\theta,\sum_{k=1}^K h_kv_k)$.
            \label{appxstmt:3}
    \end{enumerate}
    \begin{proof}
        First we establish that there exists a number $\tilde\epsilon>0$ such that, for
        each $k$,
        \begin{equation}\label{ineq:Lthetaplusevkfinite}
            L(\theta+\tilde\epsilon v_k)>-\infty.
        \end{equation}
        This fact follows from continuity of $L$ together
        with~(\ref{condition:v_k_direction_in_domain}) and the assumption that
        $L(\theta)>-\infty$. The first statement~(\ref{appxstmt:1}) in the theorem above
        then follows from Lemma~\ref{lem:Rockafellar}.
        
        Next, we establish that for any $h_1,\ldots,h_K\geq0$, there exists
        $\epsilon_h>0$ such that
        \begin{equation}\label{condition:sum_hk_vk_direction_in_domain}
            L(\theta+\epsilon_h\sum_{k=1}^{K}h_kv_k)>-\infty.
        \end{equation}
        Indeed, defining $$\epsilon_h:=\frac{\tilde\epsilon}{\sum_{k=1}^{K}h_k}$$ we see
        that
        $$\theta+\epsilon_h\sum_{k=1}^{K}h_kv_k=\theta+\frac{\sum_{k=1}^{K}h_k\tilde\epsilon
        v_k}{\sum_{k=1}^{K}h_k}$$ is a convex combination of the points
        $\theta+\tilde\epsilon v_1,\ldots,\theta+\tilde\epsilon v_k$.
        Finiteness~(\ref{condition:sum_hk_vk_direction_in_domain}) then follows from
        concavity of $L$ and finiteness~(\ref{ineq:Lthetaplusevkfinite}) of
        $L(\theta+\tilde\epsilon v_k)$ for each $k$.
        Thus, Lemma~\ref{lem:Rockafellar} gives that $DL(\theta, \sum_{k=1}^{K}h_kv_k)$
        exists and is finite, verifying statement~(\ref{appxstmt:2}) above.

        Finally, we confirm statement~(\ref{appxstmt:3}) by an application of the
        Leibnitz Integral Rule: we have equality
        \begin{align}
            \nonumber
            DL(\theta,\sum_{k=1}^{K}h_kv_k)
            &= \frac{\partial}{\partial t}\left[L\left(\theta+t\sum_{k=1}^{K}h_kv_k\right)\right|_{t=0^+} \\
            \nonumber
            &= \frac{\partial}{\partial t} \left[
                \int\left(\theta+t\sum_{k=1}^{K}h_kv_k\right)d\hat P
                - \int e^{\theta+t\sum_{k=1}^{K}h_kv_k}dM+1
                \right|_{t=0^+} \\
            \nonumber
            &= \sum_{k=1}^{K}h_k\int v_kd\hat P-\frac{\partial}{\partial t}\int
                e^{\theta+t\sum_{k=1}^{K}h_kv_k}dM
                \biggr|_{t=0^+} \\
            \label{equality:leibnitz}
            &= \sum_{k=1}^{K}h_k\int v_kd\hat P
                -\int\frac{\partial}{\partial t}
                e^{\theta+t\sum_{k=1}^{K}h_kv_k}dM
                \biggr|_{t=0^+} \\
            \nonumber
            &= \sum_{k=1}^{K}h_k\int v_kd\hat P
                -\int\sum_{k=1}^{K}h_kv_k
                e^{\theta+t\sum_{k=1}^{K}h_kv_k}dM
                \biggr|_{t=0^+} \\
            \nonumber
            &= \sum_{k=1}^{K}h_k\int v_kd\hat P
                -\sum_{k=1}^{K}h_k\int v_k
                e^{\theta}dM \\
            \nonumber
            &= \sum_{k=1}^{K}h_k 
                \frac{\partial}{\partial t} \left[
                \int\left(\theta+tv_k\right)d\hat P
                - \int e^{\theta+tv_k}dM+1
                \right|_{t=0^+} \\
            \nonumber
            &= \sum_{k=1}^{K}h_k \frac{\partial}{\partial t}\left[L\left(\theta+tv_k\right)\right|_{t=0^+}
            = \sum_{k=1}^{K}h_kDL(\theta,v_k)
        \end{align}
        where the interchange~(\ref{equality:leibnitz}) of differentiation and
        integration is allowed because, for $t\in[0,\epsilon_h]$, the function
        $e^{\theta(x)+t\sum_{k=1}^{K}h_kv_k(x)}f(x)$ is continuous and its partial
        derivative $\frac{\partial}{\partial
        t}e^{\theta(x)+t\sum_{k=1}^{K}h_kv_k(x)}f(x)=\sum_{k=1}^{K}h_kv_k(x)e^{\theta(x)+t\sum_{k=1}^{K}h_kv_k(x)}f(x)$
        is also continuous.
    \end{proof}
\end{theorem}

\nocite{*} 

\bibliography{references}{}
\bibliographystyle{plain}
\end{document}